\documentclass{amsart}
\usepackage{amsmath,amssymb,amsfonts,amsthm}
\usepackage{epsfig}
\usepackage{mathrsfs}
\usepackage{graphicx,color}

\newtheorem{theorem}{Theorem}[section]
\newtheorem{lemma}{Lemma}[section]
\theoremstyle{remark}
\newtheorem{remark}{Remark}[section]

\theoremstyle{definition}

\newtheorem{Assumption}{Assumption}

\def\R{\mathbb{R}}
\def\E{\mathbb{E}}
\def\P{\mathbb{P}}
\def\F{\mathbb{F}}

\renewcommand{\d}{\mathrm{d}}
\newcommand{\norm}[1]{\lVert#1\rVert}
\newcommand{\ind}[1]{\mathbf{1}_{#1}}
\newcommand{\abs}[1]{\left\vert#1\right\vert}
\newcommand{\ud}{\mathrm{d}}
\numberwithin{equation}{section}
\begin{document}
\title{Adapted integral representations of random variables}
\author{Georgiy Shevchenko$^*$}

\address{Department of Mechanics and Mathematics, Taras Shevchenko National University of Kyiv\\
Volodymirska 60, 01601 Kyiv, Ukraine\\
$^*$E-mail: zhora@univ.kiev.ua}

\author{Lauri Viitasaari}
\address{Department of Mathematics and System Analysis, Aalto University School of Science, Helsinki\\
P.O. Box 11100, FIN-00076 Aalto,  Finland\\
E-mail: lauri.viitasaari@aalto.fi} 

\begin{abstract}
We study integral representations of random variables with respect to general H\"older continuous processes and with respect to two particular cases; fractional Brownian motion and mixed fractional Brownian motion. We prove that arbitrary random variable can be represented as an improper integral, and that the stochastic integral can have any distribution. If in addition the random variable is a final value of an adapted H\"older continuous process, then it can be represented as a proper integral. It is also shown that in the particular case of mixed fractional Brownian motion, any adapted random variable can be represented as a proper integral.


\end{abstract}
\maketitle

\section{Introduction}

Let $(\Omega,\mathcal{F},\mathbb{F} = \{\mathcal{F}_t,t\in [0,1]\}, P)$ be a stochastic basis, and $\{X(t),t\in [0,1]\}$ be an $\mathbb{F}$-adapted process. 

We consider representations of the form 
\begin{equation}
\label{eq:intro}
\xi = \int_0^1 \psi(s) \d X(s),
\end{equation}
where $\phi$ is an $\mathbb{F}$-adapted process and $\xi$ some given $\mathcal{F}_1$-measurable random variable. 

While such representations also has theoretical interest the question is particularly motivated by mathematical finance. Indeed, $\xi$ can be viewed as the claim to be hedged and the integral representation corresponds to the value of the hedging portfolio. However, there is a constant representing the value of the claim missing in equation (\ref{eq:intro}). Consequently, claims $\xi$ with representation (\ref{eq:intro}) can be hedged with zero cost. In particular, the results presented in this paper indicates that models where the stock process $X(t)$ is H\"older continuous of some order $\alpha>\frac{1}{2}$ are rarely good models since there will be arbitrage present with relatively simple trading strategies $\psi$. 

The representations similar \eqref{eq:intro} were considered by many authors, we cite here only the most relevant results. 
The first results of this kind were established for $X=W$, the standard Wiener process. In this case, the classical It\^o representation theorem provides the representation \eqref{eq:intro} for square integrable centered random variables $\xi$ with the integrand satisfying $\int_0^1 \E \psi(t)^2 \d t<\infty$. Such representation  was shown to take place for any random variable $\xi$ in \cite{dudley}, but with integrand satisfying $\int_0^1 \psi(t)^2 \d t<\infty$ a.s. The case where $X=B^H$, a fractional Brownian motion with $H\in(1/2,1)$, was considered first in \cite{msv}. Under assumption that there exists a H\"older continuous adapted process $\{z(t),t\in[0,1]\}$ such that $z(1) = \xi$, it was shown that $\xi$ can be represented in the form \eqref{eq:intro}.  In \cite{lauri} this result was extended to a larger class of Gaussian processes, and in \cite{she-vii}, under a similar assumption, the existence of representation \eqref{eq:intro} with integrand $\psi\in C[0,1)$ was established.

In this article we generalize the results of both \cite{lauri} and \cite{she-vii} by showing the existence of the representation \eqref{eq:intro} with $\psi\in C[0,1)$ for a generic H\"older continuous process $X$ satisfying some small ball estimates.   We also show that in the case of mixed fractional Brownian motion, i.e. where $X = W+B^H$,   the representation \eqref{eq:intro} takes place for any random variable $\xi$. The structure of the article is following. Section~\ref{sec:prelim} contains basic information on the generalized Lebesgue--Stieltjes integral.  Section~\ref{sec:aux} is devoted to the auxiliary construction of processes which play an important role in proving the main representation results. Section~\ref{sec:main} contains the main results concerning the representation of random variables. 

\section{Generalized Lebesgue--Stieltjes integral}\label{sec:prelim}
This section gives a basic information on the generalized Lebesgue--Stieltjes integral, more details can be found in \cite{zahle}.
For functions $f,g\colon [a,b] \to \R$ and 
$\beta\in(0,1)$, define the fractional derivatives
\begin{gather*}
\big(D_{a+}^{\beta}f\big)(x)=\frac{1}{\Gamma(1-\beta)}\bigg(\frac{f(x)}{(x-a)^\beta}+\beta
\int_{a}^x\frac{f(x)-f(u)}{(x-u)^{\beta+1}}\d u\bigg),\\
\big(D_{b-}^{1-\beta}g\big)(x)=\frac{e^{-i\pi
\beta}}{\Gamma(\beta)}\bigg(\frac{g(x)}{(b-x)^{1-\beta}}+(1-\beta)
\int_{x}^b\frac{g(x)-g(u)}{(u-x)^{2-\beta}}\d u\bigg).
\end{gather*}
Assuming that $D_{a+}^{\beta}f\in L_1[a,b], \ D_{b-}^{1-\beta}g_{b-}\in
L_\infty[a,b]$, where $g_{b-}(x) = g(x) - g(b)$, the generalized  Lebesgue-Stieltjes integral $\int_a^bf(x)\d g(x)$ is defined as
\begin{equation*}\int_a^bf(x)\d g(x)=e^{i\pi\beta}\int_a^b\big(D_{a+}^{\beta}f\big)(x)\big(D_{b-}^{1-\beta}g_{b-}\big)(x)\d x.
\end{equation*}
From the definition, we have an immediate estimate
\begin{equation}\label{integr-estim}
\abs{\int_a^b f(x) \d g(x)}\le C\norm{f}_{\beta;[a,b]} \Lambda_\beta(g), 
\end{equation}
where 
\begin{align*}
&\norm{f}_{\beta;[a,b]} = \int_a^b \left(\frac{f(t)}{(t-a)^\beta} + \int_a^t 
\frac{\abs{f(t)-f(s)}}{(t-s)^{\beta + 1}}\d s\right)\d t,\\
&\Lambda_\beta(g) = \sup_{a\le u<v\le b} \left(\frac{\abs{g(v)-g(u)}}{(v-u)^{1-\beta}} + \int_u^v \frac{\abs{g(u)-g(z)}}{(z-u)^{2-\beta}}\d z\right).
\end{align*}
Here and in the rest of the article we will use the symbol  $C$ to denote a positive constant, whose value is of no importance and may change from one line to another.

It is easy to see that if $g$ is $\alpha$-H\"older continuous on $[a,b]$ and $\beta\in(1-\alpha,1)$, then $\Lambda_\beta(g)<\infty$. Therefore, it is possible to define $\int_a^b f(x)\d g(x)$ in the generalized Lebesgue--Stieltjes sense once the integrand $f$ satisfies $\norm{f}_{\beta;[a,b]}<\infty$. In what follows we will consider the  functions satisfying this conditions to be our admissible integrands.

\section{Adapted processes which integrate to infinity}\label{sec:aux}

At the heart of each representation lies an auxiliary construction of an adapted integrand $\psi$ such that for each $t<1$ the integral $v_t(\psi):=\int_0^t \psi(s) \d X(s)$ is finite, but the integral $\int_0^1 \psi(s)\d X(s)$ is infinite. The latter property can have different precise meanings: either $v_t(\psi) \to +\infty, t\to 1-$ or $\liminf_{t\to 1-} v_t(\psi) = -\infty$, $\limsup_{t\to 1-} v_t(\psi) = +\infty$.

\subsection{Construction in a generic case}
To obtain such auxiliary construction for general process there is essentially two key features which we study here; the process is assumed to be H\"older continuous for some order $\alpha >1/2$ and there should be some kind of estimate for small ball probability for the increment of the process. We also wish to emphasize that these properties are needed only close to the end point $t=1$ (or more generally, $t=T$). Consequently, the replication procedure can be done in arbitrary small amount of time. This can be useful for example in financial applications since one can simply wait and observe the process and study whether it might indeed have the needed properties, and then start the replication  procedure just before the ending point. For more detailed discussion in Gaussian case we refer to \cite{lauri}.

\begin{Assumption} 
\label{assu:holder}
There exist a constants $\alpha>\frac12$ such that for every $s,t\in[0,1]$ it holds
\begin{equation*}
\abs{X(t)-X(s)}\le C\abs{t-s}^{\alpha}.
\end{equation*}
\end{Assumption}

\begin{Assumption}
\label{assu:smallball}
There exists a constant $\delta>0$ such that for every $s,t\in[1-\delta,1]$ with $t=s+\Delta$ it holds
\begin{equation}
\label{smallball}
\P(\sup_{s\leq u \leq t} |X(u)- X(s)| \leq \epsilon) \leq \exp \left(-C\Delta\epsilon^{-\frac{1}{\alpha}}\right)
\end{equation}
provided that $\epsilon \leq \Delta^\alpha$.
\end{Assumption}
Note that the given upper bound for small ball probability is the usual one for many Gaussian processes and hence we wish to use this form. For example, many stationary Gaussian processes or Gaussian processes with stationary increments satisfy the given assumption. In particular, fractional Brownian motion satisfies the given assumption. For more detailed discussion on the assumption, see \cite{lauri}. We also remark that by examining our proofs below it is clear that one could relax the assumption by giving less sharp upper bound in terms of $\Delta$ and $\epsilon$ (see remark \ref{rem:smallball}).
\begin{lemma}
\label{lemma:aux_general}
Assume that the process $X$ satisfies Assumptions \ref{assu:holder} and \ref{assu:smallball}. Then there exists a $\mathbb{F}$-adapted continuous process $\phi$ on $[0,1)$ such that $\phi(0) = 0$, the integral
\begin{equation*}
\int_0^t \phi(s)\ud X(s)
\end{equation*}
exists for every $t<1$ and 
\begin{equation}
\label{rep:aux-lemma}
\lim_{t\to 1-} \int_0^t \phi(s)\ud X(s) = +\infty
\end{equation}
almost surely.
\end{lemma}
It turns out that the construction presented in the particular case of fBm in authors previous work \cite{she-vii} works for general H\"older continuous processes under our small ball assumption \ref{assu:smallball}. Hence we simply present the key points of the proof.
\begin{proof}
Fix numbers $\gamma\in\left(1,\frac{1}{\alpha}\right)$, $\eta\in\left(0,\frac{1}{\gamma \alpha}-1\right)$ and $\mu> \frac1{\alpha(1+\eta)}$. 
Set $t_0 = 0$ and $t_n = \sum_{k=1}^{n}(\Delta_{k} + \widetilde\Delta_{k})$, where $\Delta_k = K k^{-\gamma}$, $\widetilde\Delta_k = K k^{-\mu}$, $K = \big(\sum_{k=1}^\infty (k^{-\gamma}+k^{-\mu})\big)^{-1}$. Also set $t_n' = t_{n-1} + \Delta_n$, $n\ge 1$. Clearly, $t_{n-1}<t_n' < t_n$, $n\ge 1$, and $t_n\to 1$, $n\to\infty$. Note also that if $X_t$ would be $\alpha$-H\"older only close to the end point, i.e. on $[1-\delta,1]$ for some small $\delta$, then we simply set $t_1 = 1-\delta$ such that $X$ is H\"older on $[t_1,1]$ and start after $t_1$ by scaling time points properly. This also implies that the construction can be done in arbitrary small amount of time.

Next define the sequence of functions $g_n = \sqrt{x^2 + n^{-2}} - n^{-1}$, $n\ge 1$. Then $g_n(x) \uparrow |x|$, $n\to\infty$. Let also $f_n = (1+\eta) g_n(x)^\eta \frac{x}{\sqrt{x^2+n^{-2}}}$ so that $g_n(x)^{1+\eta} = \int_0^x f_n(z)\ud z$. 
For any $n\ge 1$ set
\begin{equation*}
\tau_n = \min\left\{t\geq t_{n-1} : |X(t) - X({t_{n-1}})| \geq n^{-{1}/({1+\eta})}\right\} \wedge t_n'.
\end{equation*}
Next define  
\begin{equation*}
\phi(s) = f_{n}(X(s) - X(t_{n-1}))\textbf{1}_{[t_{n-1},\tau_n)}(s)
\end{equation*}
for $s\in[t_{n-1},\tau_n]$
and 
\begin{equation*}
\phi(s) = \phi(\tau_n)\frac{\tau_n + \widetilde \Delta_n-s}{\widetilde{\Delta}_n}\textbf{1}_{(\tau_n,\tau_n+\widetilde\Delta_n]}(s)
\end{equation*}
for $s\in(\tau_n,t_{n}]$. Now by H\"older continuity of $X$ the existence of integral is clear, and we can repeat the arguments in \cite{she-vii} to obtain that
\begin{equation*}
\begin{split}
\int_0^{t_n}\phi(s)\ud X(s) &\geq 2^{-\eta}\sum_{k=1}^n |X(\tau_k)-X(t_{k-1})|^{1+\eta} \\
&- \sum_{k=1}^n k^{-1-\eta} \\
&+ \sum_{k=1}^n \int_{\tau_k}^{\tau_k+\tilde{\Delta}_k}\phi(s)\ud X(s).
\end{split}
\end{equation*}
Moreover, it is clear that the second sum converges and arguments in \cite{she-vii} imply that also the third sum converges by H\"older continuity of $X$. To conclude, the Assumption \ref{assu:smallball} implies that first sum diverges since now only finite number of events 
$$
A_n = \{\sup_{t_{n-1}\leq t\leq t'_n}|X(t) - X({t_{n-1}})| < n^{-{1}/({1+\eta})}\}
$$
happens by Borel--Cantelli Lemma and Assumption \ref{assu:smallball}. Hence the result follows.
\end{proof}
\begin{remark}
\label{rem:smallball}
By Assumption \ref{assu:smallball} we obtain that 
$$
\P(A_n) \leq \exp\left(-Cn^{\frac{1}{\alpha(1+\beta)}-\gamma}\right)
$$
for some constant $C$. Hence it is clear that our assumption on small ball probabilities could be relaxed a lot. In particular, we only need that 
$$
\sum_{n=1}^\infty \P(A_n) < \infty
$$
to apply Borel--Cantelli lemma.
\end{remark}

\subsection{Construction in pure and mixed fractional Brownian cases}
In this section we consider two important  particular cases: $X=B^H$, a fractional Brownian motion with $H>1/2$ and $X = B^H + W$, a mixed fractional Brownian motion. We start with the pure fractional Brownian case.
\begin{lemma}
\label{lemma}
Let $f(t) = (1-t)^{-H}$, $v(t) = \int_0^t f(s) dB^H(s) $. Then $\liminf_{t\to 1-} v(t) = -\infty$, $\limsup_{t\to 1-} v(t) = +\infty$ almost surely.
\end{lemma}
\begin{proof}
Define $x_n = v(1-2^{-n}) - v(1-2^{-n+1})$, $n\ge 1$. Then the sequence $\{x_n,n\ge 1\}$ is stationary Gaussian. Indeed, for any $m\ge n\ge 1$ 
\begin{align*}
\E x_n x_m
& = \alpha_H \int_{1-2^{-n+1}}^{1-2^{-n}}
\int_{1-2^{-m+1}}^{1-2^{-m}}|u-v|^{2H-2}(1-u)^{-H}(1-v)^{-H}\d u\, \d v\\
& = \alpha_H\int_{2^{-n}}^{2^{-n+1}}\int_{2^{-m}}^{2^{-m+1}}|y-x|^{2H-2}x^{-H}y^{-H}\d x\, \d y\\
& = \alpha_H\int_{1}^{2}\int_{2^{n-m}}^{2^{n-m+1}}\big|2^{-n}z- 2^{-n}w \big|^{2H-2}2^{nH}w^{-H} 2^{nH}z^{-H} 2^{-n}\d w\, 2^{-n}\d z\\
& = \alpha_H \int_{1}^{2}\int_{2^{n-m}}^{2^{n-m+1}}|z-w|^{2H-2}w^{-H}z^{-H}\d w\, \d z = r(n-m),
\end{align*}
where $\alpha_H = H(2H-1)$. 
Moreover, it is clear that $r(k) = O(2^{k(1-H)}), k\to\infty$. Therefore, defining $S_n = x_1 + x_2 + \dots + x_n$, $n\to\infty$, we have $\limsup_{n\to\infty} S_n = +\infty$, $\liminf_{n\to\infty} S_n = -\infty$ a.s. by the law of iterated logarithm for weakly dependent stationary sequences. Observing that $v(1-2^{-n}) = S_n$, we get the statement.
\end{proof}

Further we move to the case of a mixed fractional Brownian motion. This means that $X=B^H+W$, where $B^H$ is a fractional Brownian motion with $H\in(1/2,1)$, and $W$ is a standard Wiener process. Usually it is assumed that $B^H$ and $W$ are independent, but we do not impose any assumptions of such kind. Note that we understand the integral w.r.t. the standard Wiener process $W$ in the classical It\^o sense and the integral w.r.t.\ $B^H$ in the generalized Lebesgue--Stieltjes sense. 

The following lemma provides an ``auxiliary'' construction in the case despite it will not be used in the following, we give it for two reasons: to make our presentation complete and to disclose the main idea behind the proof of our main result in the mixed case. 
\begin{lemma}
Let $f(t) = (1-t)^{-1/2}$, $v(t) = \int_0^t f(s) \d(W(s) + B^H(s)) $. Then $\liminf_{t\to 1-} v(t) = -\infty$, $\limsup_{t\to 1-} v(t) = +\infty$ almost surely.
\end{lemma}
\begin{proof}
Define $u(t) = \int_0^t f(s) \d W(s)$. Then it is easy to see that $u$ has the same distribution as the time-changed Wiener process, $\{u(t), t\in [0,1]\} \overset{d}{=} \{W(-\ln(1-t)), t\in[0,1]\}$. Hence we get by the law of iterated logarithm $\liminf_{t\to 1-} u(t) = -\infty$, $\limsup_{t\to 1+} u(t) = +\infty$. So it remains to prove that the integral $\int_0^t (1-s)^{-1/2} \d B^H(s)$ is bounded. But the integrand is non-random, so the integral coincides with the so-called Wiener integral, and  its boundedness  follows from the finiteness of 
\begin{align*}
&\E \left(\int_0^1 f(s)\d B^H(s) \right)^2 \\&= H(2H-1)\int_0^1 \int_0^1 (1-t)^{-1/2}(1-s)^{-1/2}\abs{t-s}^{2H-2}\d u\,\d s\\& = 2H\,\mathrm B(2H-1,1/2).\qedhere
\end{align*}
\end{proof}

\section{Representation of random variables}\label{sec:main}

In the case of fBm it was shown in \cite{msv} that the integral $\int_0^1\phi(s)\ud B^H(s)$ can have any distribution and later in \cite{lauri} the same result was proved for wider class of Gaussian processes. Similarly, any random variable can be represented as an improper integral in these models. These results are consequence of the auxiliary construction and hence we can obtain similar results by applying auxiliary construction introduced in previous section for any H\"older process which has some small ball estimates. More precisely, a direct consequence of Lemma \ref{lemma:aux_general} is that the integral can have any distribution and if in addition we have diverging auxiliary construction on any (suitable) subinterval, then any measurable random variable can be represented as an improper integral. These results are the topic of next theorems.
\begin{theorem}
Let the process $X(t)$ satisfy Assumptions \ref{assu:holder}, \ref{assu:smallball}, and let there exist $v\in(1-\delta,1)$ such that the random variable $X(v)$ has continuous distribution. Then for any distribution function $F$ there exists a $\mathbb{F}$-adapted process $\varphi$ such that the integral
$$
\int_0^1 \varphi(s)\ud X(s)
$$
exists and has distribution $F$.
\end{theorem}
\begin{proof}
Since $X(v)$ has continuous distribution with cdf $F_X$, then $U=F_X(X(v))$ is uniformly distributed random variable and consequently, $F^{-1}(U)$ has distribution $F$. Hence it suffices to construct $\varphi$ such that
$$
\int_0^1 \varphi(s)\ud X(s) = F^{-1}[F_X(X(v))].
$$ 
Denote by $g(x) = F^{-1}[F_X(x)]$. Let $\phi$ be the process constructed in Lemma \ref{lemma:aux_general} and set $y_t = \int_v^t\phi(s)\ud X(s)$. Then $y_t \rightarrow \infty$ as $t\rightarrow 1-$. Put $\tau = \inf\{t\geq v: y_t=|g(X(v)|\}$ and
$$
\varphi(t)=\phi(t)\operatorname{sgn}g(X(v))\textbf{1}_{[v,\tau]}.
$$
Clearly $\int_0^1 \varphi(s)\ud X(s)$ has distribution $F$ and the existence of integral is obvious from which the result follows. 
\end{proof}
To replicate a distribution we needed an additional assumption that $X(v)$ has continuous distribution for some $v$. Similarly, in order to replicate arbitrary random variable we need different additional assumption. Namely, we assume that the filtration $\mathbb{F}$ is left-continuous at 1, i.e.\ $\sigma(\bigcup_{t<1} \mathcal F_t) = \mathcal F_1$.
\begin{theorem}
Assume that $\mathbb{F}$ is left-continuous at $1$ and let the process $X$ satisfy Assumptions \ref{assu:holder} and \ref{assu:smallball}. Then for any $\mathcal{F}_1$-measurable random variable $\xi$ there exists a process $\psi(s)$ such that
$$
\int_0^t \psi(s) \ud X(s)
$$
exists for every $t<1$ and
\begin{equation}
\label{rep_improper_general}
\lim_{t\rightarrow 1} \int_0^t \psi(s) \ud X(s) = \xi
\end{equation}
almost surely.
\end{theorem}
\begin{proof}
Note first that by modifying the proof of Lemma \ref{lemma:aux_general} we see that Assumption \ref{assu:smallball} implies the existence of auxiliary construction on every subinterval $[u,v]\subset[1-\delta,1]$, i.e. for every such interval there exists a process $\phi_{u,v}$ such that $\lim_{t\rightarrow v}\int_u^v \phi_{u,v}(s)\ud X(s) = \infty$. The rest follows by arguments in \cite{msv} and we only present the main steps. Define $z(t) = \tan\E[\arctan \xi |\mathcal{F}_t]$. Now by left-continuity of $\mathbb{F}$ and martingale convergence theorem we have $z(t)\rightarrow \xi,\quad t\rightarrow 1-$. Let next $t_n$ be arbitrary increasing sequence converging to $1$, and let $\phi_{t_n,t_{n+1}}$ be a process constructed in Lemma \ref{lemma:aux_general} such that $v_{t}^n = \int_{t_n}^t \phi_{t_n,t_{n+1}}(s)\ud X(s) \rightarrow \infty$ as $t\rightarrow t_{n+1}-$. Defining $\tau_n = \min\{t\ge t_n : v_t^n = |z(t_n)-z(t_{n-1})|\}$ and 
$$
\psi(s) =\sum_{n=1}^\infty \phi_{t_n,t_{n+1}}(s)\textbf{1}_{[t_n,\tau_n]}(s)\text{sign}(z(t_n)-z(t_{n-1}))
$$ 
it is clear that $\int_0^{t_n} \psi(s)\ud X(s) = z(t_{n-1})$ and on $t\in[t_n,t_{n+1}]$ the value $\int_0^t \psi(s) \ud X(s)$ is between $z(t_{n-1})$ and $z(t_n)$. Hence it follows that we have (\ref{rep_improper_general}). The existence of the integral can be shown as in the proof of Lemma \ref{lemma:aux_general}.
\end{proof}
\begin{remark}
We remark that it is also possible to construct a continuous process $\psi$ on $[0,1)$ such that $\lim_{t\rightarrow 1}\int_0^t \psi(s) \ud X(s) = \xi$ by applying similar techniques as in \cite{she-vii} or in the proof of Theorem \ref{thm:proper_general}. More precisely, after stopping $\tau_n$ let $\psi(s)$ go to zero linearly on $t\in[\tau_n,\tau_n+\Delta_n]$ for small enough $\Delta_n$, and then compensate the error arising from linear parts by setting 
$\tau_n = \min\{t\ge t_n : v_t^n = |z(t_n)-z(t_{n-1})-\int_{\tau_{n-1}}^{\tau_{n-1}+\Delta_n}\psi(s)\ud X(s)|\}$. The details are left to the reader.
\end{remark}
A particularly interesting question for us is which random variables can be represented as a proper integral.

\subsection{A proper representation in a generic case}
It turns out that with general H\"older process satisfying our small ball assumption one can represent all random variables that can be viewed as an end value of some $a$-H\"older process with arbitrary $a>0$. In the particular case of fBm this was proved first in \cite{msv}. Similar result for more general Gaussian process was derived in \cite{lauri}. However, in this case it was proved that only values $a>1-\alpha$ can be covered where $\alpha$ is the H\"older index of the process $X$. The benefit of using continuous integrands is that then one can drop unnecessary extra assumptions. Moreover, then one can cover all values of $a>0$ also in the case of general Gaussian process. More precisely, in authors previous work \cite{she-vii} it was proved that in the case of fBm one can construct a continuous integrand $\Psi$ on $[0,1)$ such that 
$$
\xi = \int_0^1 \Psi(s) \ud B^H(s).
$$ 
By examining the proof however, we obtain that only required facts are H\"older continuity, small ball estimate and the auxiliary construction with continuous integrand. Hence the arguments presented in \cite{she-vii} implies same result for our general case. Note also that, as before, the replication can be done in arbitrary small amount of time and the assumed properties are needed only close to the ending point $t=1$. 
\begin{theorem}
\label{thm:proper_general}
Let the process $X(t)$ satisfy Assumptions \ref{assu:holder} and \ref{assu:smallball}. Furthermore, assume there exists an $\mathbb{F}$-adapted process $\{z(t), t\ge 0\}$ having  H\"{o}lder continuous  paths of order
$a>0$ and such that $z(1) = \xi$. Then there exists an $\mathbb{F}$-adapted process $\{\psi(t), t\in[0,1]\}$  such that  $\psi \in C[0,1)$ a.s.\ and
\begin{equation}\label{repres}
\int_0^1 \psi(s)\ud X(s) = \xi
\end{equation}
almost surely.
\end{theorem}
The proof follows arguments presented in \cite{she-vii} but here we will give more instructive proof while some technical steps are omitted. 

The idea of the proof is to define a sequence of time points $(t_n)_{n=0}^\infty$ converging to $1$ and then track the H\"older process $z(t)$ along this sequence such that 
\begin{equation}
\label{equation_aim}
\int_0^{t_n} \psi(s)\ud X(s) = z(t_{n-1}).
\end{equation}
More precisely, we apply our diverging auxiliary construction to ``get into the right track'', and afterwards we aim to stay on this right track. Now there is two options; given that we have (\ref{equation_aim}) for some $n$ we either manage to stay on the right track and consequently we have (\ref{equation_aim}) also for $n+1$ or if we do not, then we apply the auxiliary construction together with stopping again to get ``back to the track''. Note that while we could apply the auxiliary construction separately on every interval $[t_{n-1},t_n]$, consequently the integral $\int_0^1 \psi(s)\ud X(s)$ over whole interval would not exists. Hence to obtain the result we simply have to show that we indeed manage to stay on the right path in most of the cases, and the auxiliary construction is needed only finite number of times. Finally, in order to construct a continuous integrand we simply pace to zero linearly after every time step before starting to act on the next time interval. 
\begin{proof}[Proof of Theorem \ref{thm:proper_general}]
Choose some $\beta\in(1-\alpha,1)$. Let $\Delta_k$ be sequence such that $\sum_{k=1}^\infty \Delta_k=1$ and define time points $t_0=0$, $t_n = \sum_{k=1}^n \Delta_k$. Set also $t'_n = t_{n-1} + \frac{\Delta_n}{2}$. Note that now $t_{n-1}<t'_n<t_{n}$. Let also $\widetilde{\Delta}_k$ be a sequence to be determined later such that $\widetilde\Delta_k \leq \frac{\Delta_k}{2}$. Following the idea described above, our aim is to define continuous integrand such that we track the process $z(t)$ on intervals $[t_{n-1},t'_n]$ and then we go linearly to zero such that the integrand hits zero before time $t_{n+1}$. Then on $[t_{n+1},t'_{n+1}]$ we define continuous integrand $\psi$ such that $\psi(t_{n+1})=0$ and we are tracking the process $z(t)$. Note also that on every step we have to compensate the error arising from linear parts. We will first explain naively the construction which is somewhat simple. The end of the proof is devoted to analysis on different parameters where we show that one can indeed chose them such that we obtain our result.\\
\textbf{Step 1. Construction}.  
We start by setting $\psi(t) = 0$ on the interval $[t_0,t_1]$. Moreover, we set $\tau_1 = t_1$. 

Denote $y(t) = \int_0^t \psi(s) \ud X(s)$, $\xi_n = z(t_{n-1})$, $g_{n}(x) = \sqrt{x^2+\epsilon_n^2}-\epsilon_n$ for some sequence $\epsilon_n$ and let now $n\ge 2$. 
To describe our construction mathematically, we want to define the process $\psi$ on $[t_{n-1}, t_{n}]$ such that
\begin{enumerate}
\item $\psi(t_{n-1}) = \psi(t_{n}) = 0$;
\item $y(\tau_{n}) = \xi_n$ for some $\tau_{n}\in [t_{n},t'_{n}]$;
\item $\psi$ is linear on $[\tau_{n},\tau_{n}+\widetilde\Delta_n]$ and zero afterwards, i.e.
\begin{equation}\label{psilinear}
\psi(t) = \psi(\tau_n)\frac{\tau_n + \widetilde \Delta_n-t}{\widetilde{\Delta}_n}\textbf{1}_{(\tau_n,\tau_n+\widetilde\Delta_n]}(t),\quad t\in [\tau_n, t_{n}].
\end{equation}
\end{enumerate}
Now the construction is different whether we are already ``on the right path'' (case A) or not (case B) in which case we apply the auxiliary construction of Lemma \ref{lemma:aux_general}.

Case A) $y({\tau_{n-1}})=\xi_{n-1}$. For a sequence $a_n$ to be determined later, define 
\begin{equation*}
\tau_{n} = \inf\left\{t\geq t_{n-1} : a_ng_{n}(X(t)- X(t_{n-1}))= \left|\Lambda_n\right|\right\} \wedge t_{n}',
\end{equation*}
where 
$
\Lambda_n = \xi_n - y(t_{n-1}) = \xi_n - \xi_{n-1} -\int_{\tau_{n-1}}^{t_{n-1}}\psi(s)\ud X(s).
$
Put $$\psi(t) = a_n g'_{n} (X(t) - X(t_{n-1}))\operatorname{sign} \Lambda_n,\quad t\in [t_{n-1}, \tau_{n}]$$
and define it by \eqref{psilinear} on $[\tau_{n},t_{n}]$.
Now since $X(t)$ is H\"older continuous of order $\alpha>\frac{1}{2}$, it obeys the classical change of variable rule. Hence we get
$$
y(t) = y(t_{n-1}) + a_n g_{n}(X(t) - X(t_{n-1}))\operatorname{sign} \Lambda_n,\quad t\in [t_{n-1}, \tau_{n}];
$$
in particular,  $y(\tau_{n}) = \xi_n$ provided that $\tau_{n}<t'_{n}$. 

Case B)  $y({\tau_{n-1}}) \neq \xi_{n-1}$. Since Assumption \ref{assu:smallball} implies that there exists diverging auxiliary construction also on every subinterval, there exists an adapted continuous process $\{\phi_n(t),t\in[t_{n-1},t_n']\}$ such that $v_n(t):=\int_{t_{n}}^t \phi_n(s) \ud X(s)\to \infty$, $t\to t'_n-$. Therefore we can define the stopping time $\tau_{n} = \inf \{t\in [t_{n-1},t_{n}'): v(t) = \abs{ \xi_n-y(\tau_{n-1})}\}$. Then we put $\psi (t) = \phi_n(t)\operatorname{sign}(\xi_n - y(\tau_{n-1}))$, $t\in [t_{n-1},\tau_n]$, and use \eqref{psilinear} on $[\tau_{n},t_{n}]$. Clearly, $y(\tau_n) = \xi_n$.\\
\textbf{Step 2. ``Staying on the right path'' and continuity of the integral}.
To obtain our result we wish to apply Assumption \ref{assu:smallball} to obtain that we have Case A) in most of the cases and that representation (\ref{repres}) holds. For the latter one, it is sufficient to prove that the integral $\int_0^t \psi(s)\ud X(s)$ is continuous at $t=1$ which also implies the existence of the integral. Consequently, we end up to some restrictions on free parameters. Similarly, note that in order to prove that we have Case A) in most of the cases we have to prove that the event
\begin{equation*}
A_n = \Big\{\sup_{t\in[t_{n-1},t_n')} a_ng_{n}(X(t)- X(t_{n-1}))\leq |\Lambda_n|\Big\}
\end{equation*}
happens only finite number of times. Now by following arguments in \cite{she-vii} we obtain $\abs{\psi(\tau_{n-1})}\le a_n$ for every $n$ and 
\begin{align*}
\abs{\int_{\tau_{n-1}}^{t_{n-1}}\psi(s)\ud X(s)}\le C_\epsilon (\omega)\abs{\phi(\tau_{n-1})} \widetilde{\Delta}_{n-1}^{\alpha-\epsilon}\le C_\epsilon (\omega)a_n\widetilde{\Delta}_{n-1}^{\alpha-\epsilon}
\end{align*}
by H\"older continuity of $X$. Moreover, observing that $g_{\epsilon_n}(x) \geq |x|-\epsilon_n$ we obtain that the event $A_n$ implies
\begin{equation*}
\sup_{t\in[t_{n-1},t_n')} |X(t)- X(t_{n-1}| \leq a_n^{-1}|\xi_n-\xi_{n-1}| + C_{\epsilon}(\omega)\widetilde\Delta_{n-1}^{\alpha-\epsilon}+\epsilon_n.
\end{equation*}
Moreover, by H\"older continuity of $z(t)$ this implies that also
\begin{equation}
\label{small_needed}
\sup_{t\in[t_{n-1},t_n')} |X(t)- X(t_{n-1})| \leq C(\omega)a_n^{-1}\Delta_n + C_{\epsilon}(\omega)\widetilde\Delta_{n-1}^{\alpha-\epsilon}+\epsilon_n.
\end{equation}
Now the idea is to choose parameters such that (\ref{small_needed}) takes place only finite number of times. Next we will study the continuity of the integral. For this it suffices to show that 
$\int_{\tau_n}^{1} \psi(s)\ud X(s) \to 0$, $n\to\infty$, which would follow from $\norm{\psi}_{\beta, [\tau_n,1]}\to 0$, $n\to\infty$. Assume now that we have chosen parameters such that (\ref{small_needed}) takes place only finite number of times. We write
$$
\norm{\psi}_{\beta; [\tau_n,1]} = I_1 + I_2, 
$$
where 
$$
I_1 =  \int_{\tau_n}^1\frac{\abs{\psi(t)}}{(t-\tau_n)^\beta} \ud s,\quad I_2  = \int_{\tau_n}^1\int_{\tau_n}^{t}\frac{\abs{\psi(t)-\psi(s)}}{(t-s)^{\beta+1}}\ud s\,\ud t.
$$
We follow arguments presented in \cite{she-vii} to obtain bounds for terms $I_1$ and $I_2$ with our general parameters, and some technical details will be omitted.
First we estimate
\begin{align*}
I_1 & = \int_{\tau_n}^{t_n}\frac{\abs{\psi(t)}}{(t-\tau_n)^\beta}\ud t +  \sum_{k=n}^{\infty} \int_{t_k}^{t_{k+1}}\frac{\abs{\psi(t)}}{(t-\tau_n)^{\beta}}\ud t\\
& \le C a_n {\Delta}_n^{1-\beta} + C\sum_{k=n}^{\infty} \frac{a_k\Delta_k}{(t_k - \tau_{n-1})^{\beta}} \le C \sum_{k=n}^{\infty} a_k\Delta_k^{1-\beta}.
\end{align*}
For $I_2$ we write
\begin{align*}
I_2 & = 
\sum_{k=n}^\infty \int_{t_k}^{\tau_{k+1}}  \int_{\tau_n}^{t} \psi(t,s)\ud s \,\ud t 
+ \sum_{k=n}^\infty \int_{\tau_{k}}^{t_{k}}  \int_{\tau_n}^{t} \psi(t,s)\ud s \,\ud t\\
&  = \sum_{k=n}^\infty \int_{t_k}^{\tau_{k+1}}  \int_{\tau_n}^{t_k} \psi(t,s)\ud s \,\ud t 
+ \sum_{k=n}^\infty \int_{\tau_{k}}^{t_{k}}  \int_{\tau_n}^{\tau_k} \psi(t,s)\ud s \,\ud t\\
& + \sum_{k=n}^\infty \int_{t_k}^{\tau_{k+1}}  \int_{t_k}^{t} \psi(t,s)\ud s \,\ud t 
+ \sum_{k=n}^\infty \int_{\tau_{k}}^{t_{k}}  \int_{\tau_k}^{t} \psi(t,s)\ud s \,\ud t\\
&=: J_1 + J_2 + J_3 + J_4,
\end{align*}
where $\psi(t,s) = \abs{\psi(t)-\psi(s)}(t-s)^{-\beta-1}$. Now arguments in \cite{she-vii} imply that we have
$$
J_1 \leq C\sum_{k=n}^{\infty} a_k\Delta_k^{1-\beta},
$$
$$
J_2 \leq C\sum_{k=n}^{\infty} a_k\Delta_k^{1-\beta},
$$
and
$$
J_4 \leq C\sum_{k=n}^{\infty} \widetilde\Delta_k^{-1}\Delta_k^{2-\beta}.
$$
Moreover, for $J_3$ we get by H\"older continuity of $X$ that
\begin{equation*}
J_3 \leq C(\omega)\sum_{k=n}^{\infty}a_k\epsilon_k^{-1}\Delta_k^{1+\alpha-\epsilon-\beta}.
\end{equation*}
To summarize, we need to choose parameters such that (\ref{small_needed}) happens only finite number of times and \begin{enumerate}
\item
$$
\sum_{k=n}^{\infty}a_k\Delta_k^{1-\beta}\rightarrow 0,
$$
\item
$$
\sum_{k=n}^{\infty}\widetilde\Delta_k^{-1}\Delta_k^{2-\beta}\rightarrow 0,
$$
\item
$$
\sum_{k=n}^{\infty}a_k\epsilon_k^{-1}\Delta_k^{1+\alpha-\epsilon-\beta}\rightarrow 0.
$$
\end{enumerate}
\textbf{Step 3. Analysis of the parameters}.
Next we prove that we can choose parameters such that we obtain $(1)-(3)$ and (\ref{small_needed}) happens only finite number of times. For simplicity let us first put $a_n=\Delta_n^{-\mu}$, $\widetilde{\Delta}_n = \Delta_n^\gamma$, and $\epsilon_n = \Delta_n^\kappa$ for some parameters $\mu$, $\gamma$ and $\kappa$. 
With these choices (\ref{small_needed}) implies that
$$
\sup_{t\in[t_{n-1},t_n')} |X(t)- X(t_{n-1})| \leq C(\omega)\Delta_n^{\lambda},
$$
where 
$$
\lambda = \min(\mu + a, \gamma(\alpha-\epsilon),\kappa).
$$
Moreover, for any small number $\hat\epsilon$ there exists $N(\omega)$ such that
$$
C(\omega)\Delta_n^{\hat\epsilon} \leq 1,\quad n\geq N(\omega).
$$
Note now that the restriction $\epsilon \leq T^\alpha$ in Assumption \ref{assu:smallball} implies that we have to choose parameters such that
$\Delta_n^{\lambda -\hat\epsilon} < \Delta_n^\alpha$, or equivalently $\lambda-\hat\epsilon> \alpha$. With such choices and applying Assumption \ref{assu:smallball} we obtain 
$$
\mathbb{P}\left(\sup_{t\in[t_{n-1},t_n')} |X(t)- X(t_{n-1})| \leq \Delta_n^{\lambda-\hat\epsilon}\right) \leq exp\left(-C\Delta_n^{1-\frac{\lambda-\hat\epsilon}{\alpha}}\right)
$$
which is clearly summable provided $\Delta_n$ converges to zero fast enough. 

Consider next restrictions $(1)-(3)$. With our choices we demand that
$$
\sum_{k=n}^{\infty}\Delta_k^{1-\beta-\mu}\rightarrow 0,
$$
$$
\sum_{k=n}^{\infty}\Delta_k^{2-\beta-\gamma}\rightarrow 0,
$$
and
$$
\sum_{k=n}^{\infty}\Delta_k^{1+\alpha-\epsilon-\beta-\mu-\kappa}\rightarrow 0.
$$
Again, these conditions are clearly satisfied provided that all the exponents are positive and $\Delta_k$ decays to zero fast enough. Hence, by collecting all restrictions, we obtain that we have to choose $\beta$, $\mu$, $\kappa$, $\gamma$, $\epsilon$, and $\hat\epsilon$ such that;
\begin{enumerate}
\item
$\mu + a - \hat\epsilon > \alpha$,
\item
$\kappa-\hat\epsilon > \alpha$,
\item
$\gamma(\alpha-\epsilon)-\hat\epsilon > \alpha$,
\item
$1-\beta-\mu>0$,
\item
$2-\beta-\kappa>0$,
\item
$1+\alpha-\epsilon-\beta-\mu-\kappa > 0$.
\end{enumerate}
The first three restrictions arises from Assumption \ref{assu:smallball} and the latter three from $(1)-(3)$. Note also that by choosing $\epsilon$ and $\hat\epsilon$ small enough we can actually omit them on the conditions $(1)-(6)$. Now $(3)$ implies that $\gamma>1$ which is consistent with $\widetilde{\Delta}_n \leq \Delta_n/2$ for $n$ large enough. Next combining $(4)$ and $(1)$ we obtain that
$$
\alpha - a < \mu < 1-\beta
$$
which is possible if we choose $\beta\in(1-\alpha,1-\alpha+a)$. Next combining $(2)$ and $(5)$ we obtain
$$
\alpha < \kappa < 2-\beta,
$$
and together with restrictions on $\beta$ this is possible if
$$
\alpha < \kappa < 1+\alpha - a.
$$
This is clearly possible since we can, without loss of generality, assume that $a<\alpha$. It remains to study $(6)$. By choosing $\mu= \alpha - a + \delta$ and $\kappa= \alpha+\delta$ for $\delta$ small enough we obtain
\begin{equation*}
\begin{split}
1+\alpha -\beta - \mu -\kappa &= 1+ \alpha - \beta -\alpha +a -\delta - \alpha - \delta \\
&= 1-\beta + a -\alpha - 2\delta\\
&> 0 
\end{split}
\end{equation*}
since $\beta < 1-\alpha + a$. To conclude, we obtained that we can choose parameters properly and it remains to choose $\Delta_n$ such that it converges to zero fast enough.
\end{proof}
\begin{remark}
In \cite{she-vii} the authors defined $\Delta_n = n^{-\nu}$ and then chose $\nu$ properly to obtain the result. Now we obtained that one can choose $\Delta_n$ in many different ways. For example, one can choose $\Delta_n = \frac{K}{2^n}$ with $K=\left(\sum_{k=1}^\infty 2^{-k}\right)^{-1}$.
\end{remark}

\subsection{Representation in the mixed case}
Next we turn to the representation w.r.t. the mixed fractional Brownian motion $B^H + W$ with $H\in(1/2,1)$. We recall that the integral w.r.t. $W$ is understood in the It\^o sense, that w.r.t.\ $B^H$, in the generalized Lebesgue--Stieltjes sense. 
\begin{theorem}
Assume that the filtration $\F$ is left-continuous at $1$. Then for  any $\mathcal F_1$-measurable random variable $\xi$ there exists an $\mathbb{F}$-adapted process $\psi$ such that 
\begin{equation}\label{mixed-representation}
\int_0^1 \psi(t) \d (B^H(t) + W(t)) = \xi 
\end{equation}
a.s.
\end{theorem}
\begin{proof}
The proof is similar to that of the main result in \cite{dudley}. Choose some $\beta\in(1-H,1/2)$. 
In view of the left-continuity of $\F$ at 1, for each $k$ there exists an $\mathcal F_{1-2^{-k}}$-measurable random variable $\xi_k$ such that $\xi_k \to \xi$, $k\to\infty$, a.s. Take a subsequence $k(n)\to \infty$, $n\to\infty$ such that 
\begin{equation}\label{n3n2}
\P(\abs{\xi_{k(n)}-\xi}>n^{-3})\le n^{-2}
\end{equation}
and denote $t_n = 1-2^{-k(n)}$.
The integrand $\psi$ will be of the form $\psi = \sum_{n=1}^\infty  (t_{n+1}-t)^{-1/2}\ind{[t_n,\tau_n)}$ with some 
$\tau_n \in(t_n, t_{n+1}]$. First we make some a priori estimates concerning the integrand. To this end, consider 
\begin{gather*}
\norm{\psi}_{\beta;[t_n,t_{n+1}]} = I_1 + I_2, 
\end{gather*}
where 
\begin{align*}
I_1 &= \int_{t_n}^{t_{n+1}}\frac{\abs{\psi(t)}}{(t-t_n)^{\beta}}\d t = 
\int_0^{\tau_n} \frac{(t_{n+1}-t)^{-1/2}}{(t-t_n)^{\beta}}\d t\le C (t_{n+1} - t_n)^{1/2-\beta};\\
I_2 &= \int_{t_n}^{t_{n+1}}\int_{t_n}^{t}
\frac{\abs{\psi(t)-\psi(s)}}{(t-s)^{\beta+1}}\d s\, \d t  \\&=  \int_{t_n}^{\tau_n} \int_{t_n}^{t}\frac{\abs{(t_{n+1}-t)^{-1/2}-(t_{n+1}-s)^{-1/2}}}{(t-s)^{\beta+1}}\d s\, \d t\\
&  + \int_{\tau_n}^{t_{n+1}}\int_0^{\tau_n}\frac{(t_{n+1}-s)^{-1/2}}{(t-s)^{\beta+1}}\d s\, \d t \\
& \le \int_{t_n}^{\tau_n} \int_{t_n}^{t}\frac{\d s \, \d t}{(t_{n+1}-t)^{1/2}(t_{n+1}-s)^{1/2}(t-s)^{\beta+1/2}} \\
& + C \int_0^{\tau_n} (t-\tau_n)^{-\beta-1/2}\d t \le C (t_{n+1}-t_n)^{1/2-\beta}.
\end{align*}
In particular, denoting $v^H_n = \int_{t_n}^{t_{n+1}} \psi(t) \d B^H(t)$, we have
\begin{equation}\label{vhbestim}
\abs{v^H_n}\le C(t_{n+1}-t_n)^{1/2-\beta}\Lambda(B^H),
\end{equation}
 whence 
\begin{align*}
\P (\abs{v_n^H}\ge n^{-3}) &\le \P\left(\abs{\Lambda(B^H)}\ge Cn^{-3}(t_{n+1}-t_n)^{\beta-1/2}\right)\\
&\le \P\left(\abs{\Lambda(B^H)}\ge Cn^{-3}2^{(1/2-\beta)n}\right)\\
&\le C \E\abs{\Lambda(B^H)} n^{3} 2^{(\beta-1/2)n},
\end{align*}
consequently, 
\begin{equation}\label{Pvn>n3}
\sum_{n=1}^\infty\P (\abs{v_n^H}\ge n^{-3})<\infty.
\end{equation}

Now we define $\tau_n$ consecutively.  Denote $v(t) = \int_0^t \psi(s) \d (B^H(s)+W(s))$. For given $n\ge 1$ assume that $\tau_k$ is defined for $k=1,\dots, n-1$ (for $n=0$, assume nothing) and denote $v_n^W(t) = \int_{t_n}^{t}(t_{n+1}-t)^{-1/2}\d W(t)$. Since $\int_{t_n}^{t_{n+1}}(t_{n+1}-t)\d t=+\infty$, we have $\liminf\int_{t_n}^{t}(t_{n+1}-t)^{-1/2}\d W(s) = -\infty$ a.s., $\limsup\int_{t_n}^{t}(t_{n+1}-t)^{-1/2}\d W(s) = +\infty$ a.s. Therefore, the stopping time $\tau_n = \inf\{t\ge  t_n :  v_n^W(t)  = \xi_{k(n)}-v(t_n)\}\wedge t_{n+1}$ satisfies $\tau_n<t_{n+1}$ a.s.

We need to show \eqref{mixed-representation}. First we argue that $\int_0^1 \psi(t)^2 \d t<\infty$ so that the integral $\int_0^1 \psi(s) \d W(s)$ is well defined. Denote $G_n = \int_{t_n}^{t_{n+1}} \psi(t)^2 \d t$. As in \cite{dudley}, we have $\P(G_n \ge n^{-2}\mid \mathcal F_{t_n}) \le n \abs{\xi_{k(n)}-v(t_n)}$. For $n\ge 2$, note that 
\begin{equation}\label{vnt}
\begin{aligned}
v(t_n)& = v(t_{n-1}) + \int_{t_{n-1}}^{t_n} \psi(t)\d W(t) + \int_{t_{n-1}}^{t_n} \psi(t)\d B^H(t)\\& = \xi_{k(n-1)} + \int_{t_{n-1}}^{t_n} \psi(t)\d B^H(t) = \xi_{k(n-1)} + v_n^H
\end{aligned} 
\end{equation}
and estimate 
\begin{align*}
\abs{\xi_{k(n)}-v(t_n)}\le \abs{\xi_{k(n)}-\xi}+\abs{\xi_{k(n-1)}-\xi} + \abs{v_n^H}.
\end{align*}
Therefore, taking into account \eqref{n3n2} and \eqref{Pvn>n3}, we get that 
\begin{align*}
&\sum_{n=1}^\infty \P(G_n \ge n^{-2})\le \sum_{n=1}^{\infty }\left( 3n^{-2} + \P(n\abs{\xi_{k(n)}-v(t_n)}\ge 3n^{-2})\right)\\
& \le \sum_{n=1}^{\infty }\bigg( 3n^{-2} + \P(\abs{\xi_{k(n)}-\xi}\ge n^{-3})\\&\qquad + \P(\abs{\xi_{k(n-1)}-\xi}\ge n^{-3})  + \P\left(\abs{v_n^H}\ge n^{-3}\right) \bigg)<\infty.
\end{align*}
Then the Borel--Cantelli lemma implies that $\int_0^1 \psi(t)^2 \d t = \sum_{n=1}^\infty G_n<\infty$ a.s.

To show the existence of the integral $\int_0^1 \psi(t) \d B^H(t)$, we write $\norm{\psi}_{\beta;[0,1]}\le \sum_{n=1}^\infty\norm{\psi}_{\beta;[t_n,t_{n+1}]}$  and use the above estimates for $\norm{\psi}_{\beta;[t_n,t_{n+1}]}$.

In view of \eqref{vhbestim} and \eqref{vnt}, 
$v(t_n)\to \xi$, $n\to\infty$. It remains to prove that $v(1)-v(t_n) = \int_{t_n}^1 \psi(t)\d (B^H(t)+W(s))\to 0$, $n\to\infty$. 
The convergence $\int_{t_n}^1 \psi(t)\d W(t)\to 0$, $n\to \infty$ follows from the convergence of the series $\sum_{n=1}^\infty G_n$, and 
$\abs{\int_{t_n}^1 \psi(t)dB^H(t)}\le \sum_{k=n}^\infty \abs{v_n^H}<\infty$ thanks to \eqref{vhbestim}, which concludes the proof.
\end{proof}
\begin{remark}
It is straightforward to generalize the statement to the case where $X(t) = \int_0^t \sigma(s) \d W(s) + Z(t)$, where the process $\sigma$ is an $\F$-adapted bounded non-vanishing process, and $Z$ is H\"older continuous of some order $\alpha>1/2$.
\end{remark}

\section*{Acknowledgments}

The first author would like to thank Dr.\ Maria Victoria Carpio-Bernido and Dr.\ Christopher Bernido for their great hospitality during his visit to Jagna in January 2014.

\bibliographystyle{ws-procs9x6}
\bibliography{compilation}

\begin{thebibliography}{1}

\bibitem{dudley}
R.~M. Dudley, {Wiener functionals as Ito integrals.}, {\em Ann. Probab.} {\bf
  5}, 140  (1977).

\bibitem{msv}
Y.~Mishura, G.~Shevchenko and E.~Valkeila, Random variables as pathwise
  integrals with respect to fractional {B}rownian motion, {\em Stochastic
  Process. Appl.} {\bf 123}, 2353  (2013).

\bibitem{lauri}
L.~Viitasaari, Integral representation of random variables with respect to
  gaussian processes  (2013), arXiv:math.PR/1307.7559.

\bibitem{she-vii}
G.~Shevchenko and L.~Viitasaari, Integral representation with adapted
  continuous integrand with respect to fractional {B}rownian motion  (2014),
  arXiv:math.PR/1403.2066.

\bibitem{zahle}
M.~Z{\"a}hle, {Integration with respect to fractal functions and stochastic
  calculus. {I}.}, {\em Probab. Theory Relat. Fields} {\bf 111}, 333  (1998).

\end{thebibliography}

\end{document}